\documentclass{llncs}

\usepackage{footnote}
\makesavenoteenv{tabular}
\makesavenoteenv{table}

\usepackage{graphicx}
\usepackage{eurosym}
\usepackage{cleveref}
\usepackage{url}

\begin{document}

\title{A compilation of LEGO Technic parts to support learning experiments on linkages}
\author{Zolt\'an Kov\'acs\inst{1} \and Benedek Kov\'acs\inst{2}}
\institute{
The Private University College of Education of the Diocese of Linz\\
Salesianumweg 3, A-4020 Linz, Austria\\
\email{zoltan@geogebra.org} \and
Europagymnasium Auhof\\
Aubrunnerweg 4, A-4040 Linz, Austria\\
\email{kovben2004@gmail.com}
}

\maketitle

\begin{abstract}

We present a compilation of LEGO Technic parts to provide easy-to-build constructions
of basic planar linkages. Some technical issues and their possible solutions are discussed.
To solve questions on fine details---like deciding whether the motion is an exactly straight line or not---we
refer to the dynamic mathematics software tool GeoGebra.

\keywords{linkages, LEGO, STEM education, GeoGebra, algebraic geometry}

\end{abstract}

\section{Making mathematics visible via linkages}
Making mathematics visible is a very recent topic of today's mathematics education.
Connecting formulas and real life problems is one of the main interests of the STEM movement.
Computer aided visualization is one possible way to support understanding, but, on
the other hand, utilizing palpable objects in education promotes better connection
between mathematical ideas and real life objects in today's electronic world, too.

One possible field of such a connection is modeling planar linkages with LEGO Technic
components. One obvious benefit of using LEGO parts is that many children already have
such a construction kit at home, and the others can buy or borrow them with the help
of their school without difficulties. On the other hand, planar linkages support
various levels of mathematical contents, including the concept of a circle,
a geometric locus, the algebraic equation of a curve, manipulation on equation systems,
factorization of polynomials, and automated theorem proving. Examples on supporting
these kinds of contents can already be found in the literature \cite{oldenburg,CGTA2017,GGG-linkages}.

In our approach the motion of the moving parts of a linkage can be visualized with the help of a sheet of
paper and a pen refill. This kind of extension of the traditional way of using LEGO
components opens new possibilities to have a closer look on the movements. On the other
hand, such constructions still cannot be precise enough to find inaccuracies in the
mathematical concept of the motion being built. To completely solve such issues
a computer should be involved---one that is fast enough to manipulate on tens of
equations in many variables.

This paper consists of 4 sections. In Section \ref{sec:overview} we consider all linkages
that can be built from our compilation, in order to introduce them as a sequence of
materialized objects for various mathematical concepts. In Section \ref{sec:build} we
give a list of possible ways to assemble the linkages, and we report on some technical
difficulties and their possible solutions. Finally, in Section \ref{sec:summary} we
discuss some thoughts on possible further extensions.

\section{An overview on the linkages in the compilation}
\label{sec:overview}

\begin{figure}
\begin{center}
\includegraphics[scale=0.3]{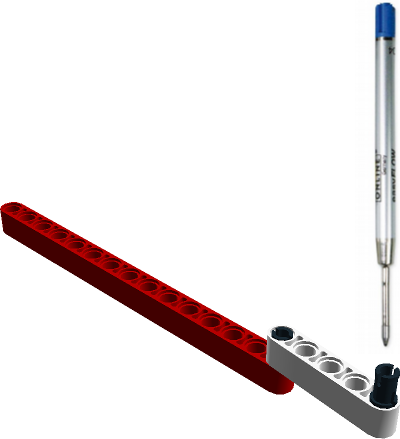}\hskip0.2cm
\includegraphics[scale=0.5]{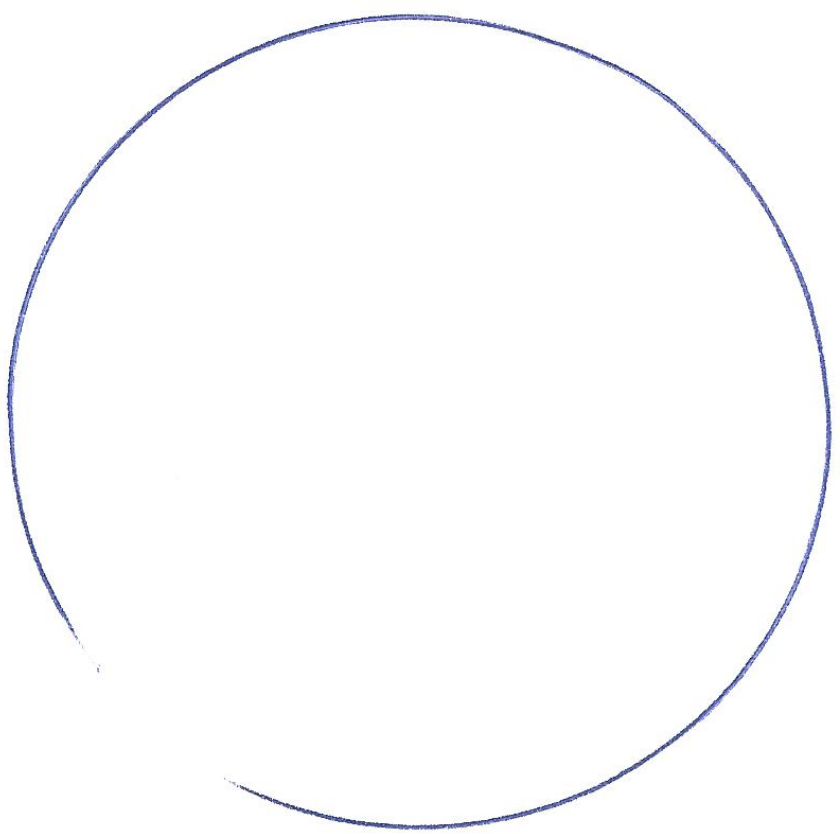}
\caption{A simple ``pair of compasses'' and the circle it draws.}
\label{fig:circle}
\end{center}
\end{figure}

\begin{figure}
\begin{center}
\includegraphics[scale=0.15]{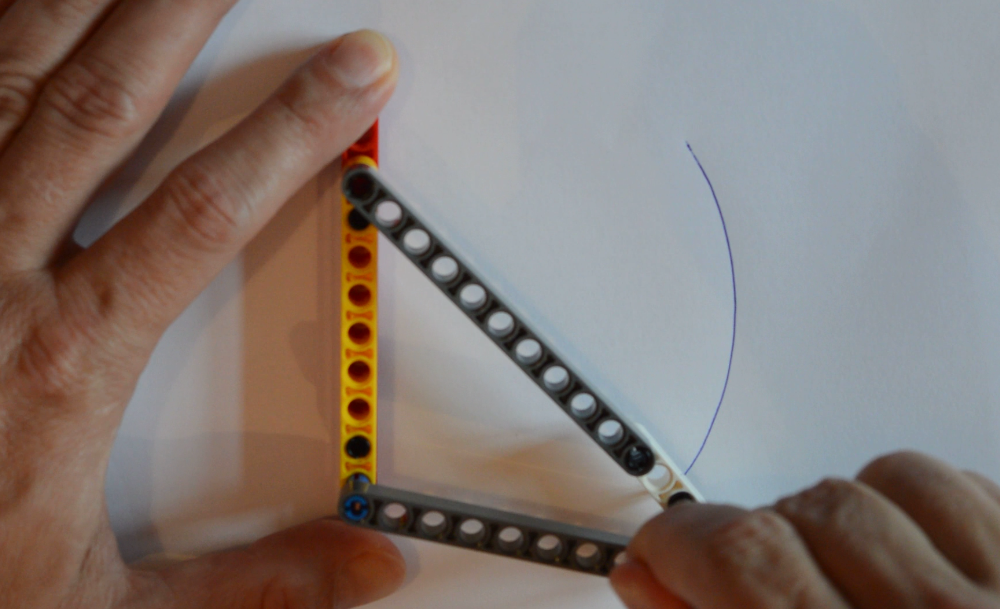}\hskip0.2cm
\includegraphics[scale=0.15]{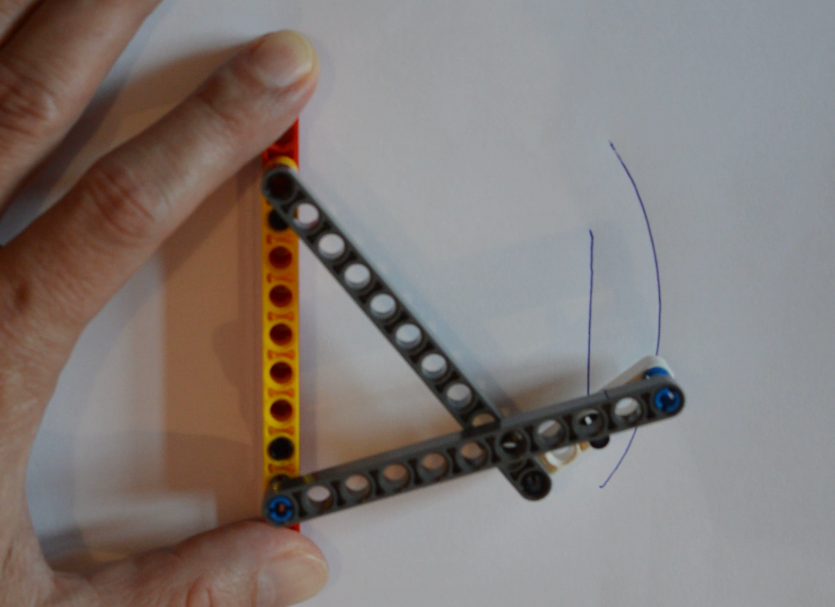}
\caption{Chebyshev's linkage built as a LEGO Technic model, and the curves it can draw.}
\label{fig:how}
\end{center}
\end{figure}

\begin{figure}
\begin{center}
\includegraphics[scale=1]{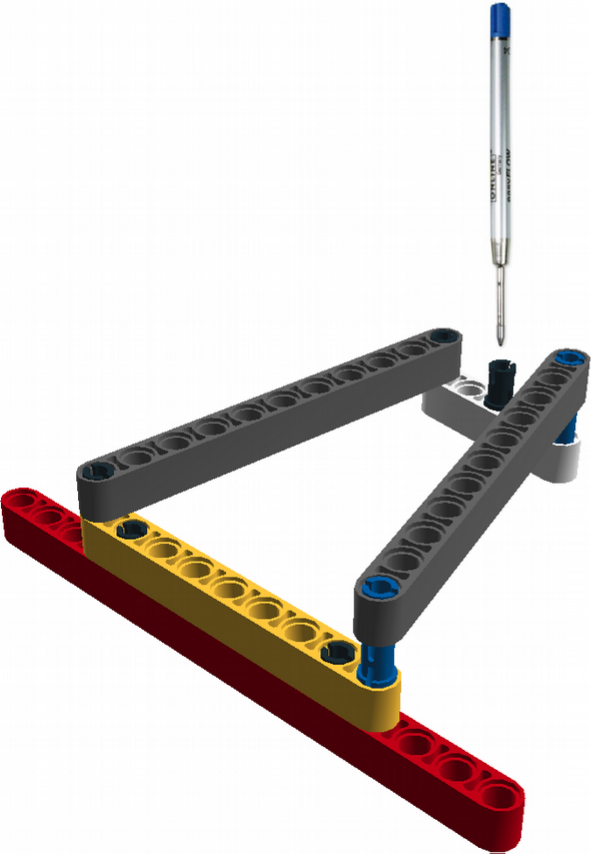}\hskip0.2cm
\raisebox{3cm}{\includegraphics[scale=0.5]{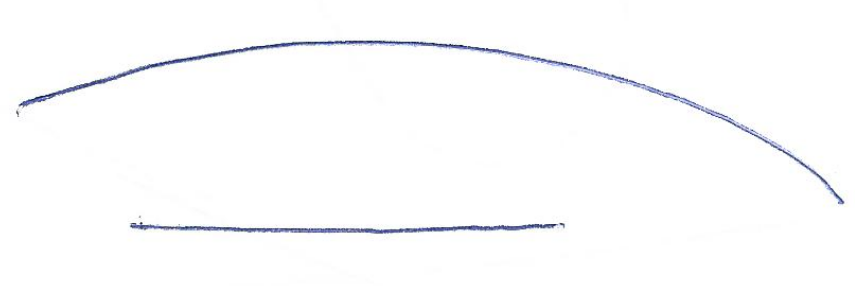}}
\caption{Chebyshev's linkage (1854) and the pieces of curves it draws.}
\label{fig:chebyshev}
\end{center}
\end{figure}

\begin{figure}
\begin{center}
\includegraphics[scale=0.5]{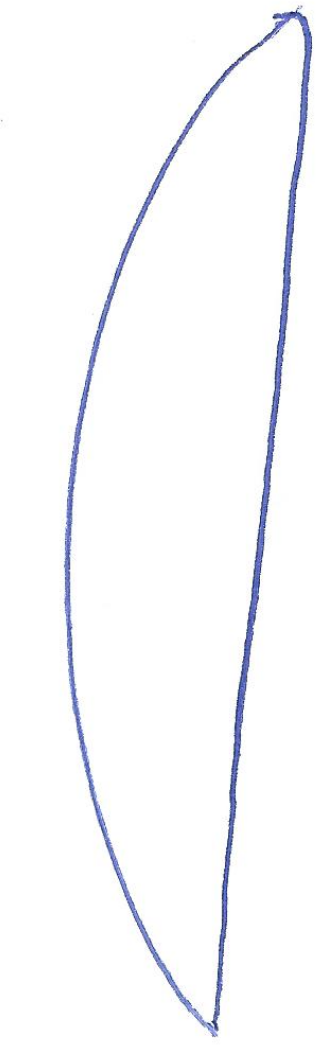}\hskip0.2cm
\includegraphics[scale=0.2]{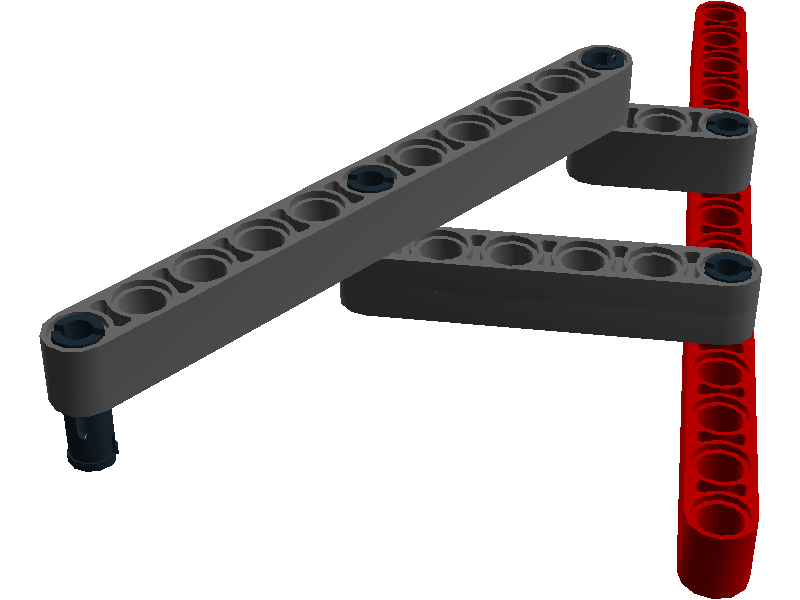}
\caption{Chebyshev's $\lambda$ mechanism (1878) and the sextic it draws in one movement.
The pen refill should be inserted at the left end of the beam with 11 holes. Note that
two pieces of half-sized beam 6 are used.}
\label{fig:chebyshevl}
\end{center}
\end{figure}

\begin{figure}
\begin{center}
\includegraphics[scale=0.3]{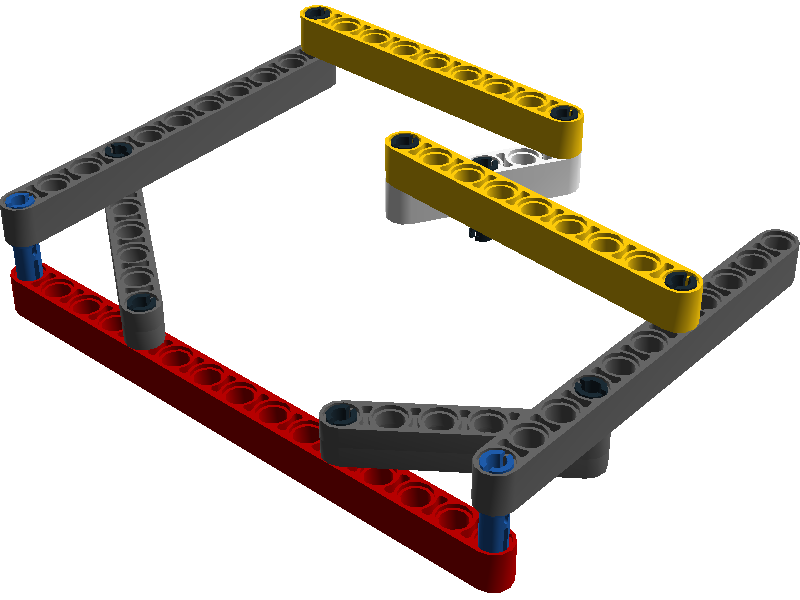}\hskip0.2cm
\includegraphics[scale=0.5]{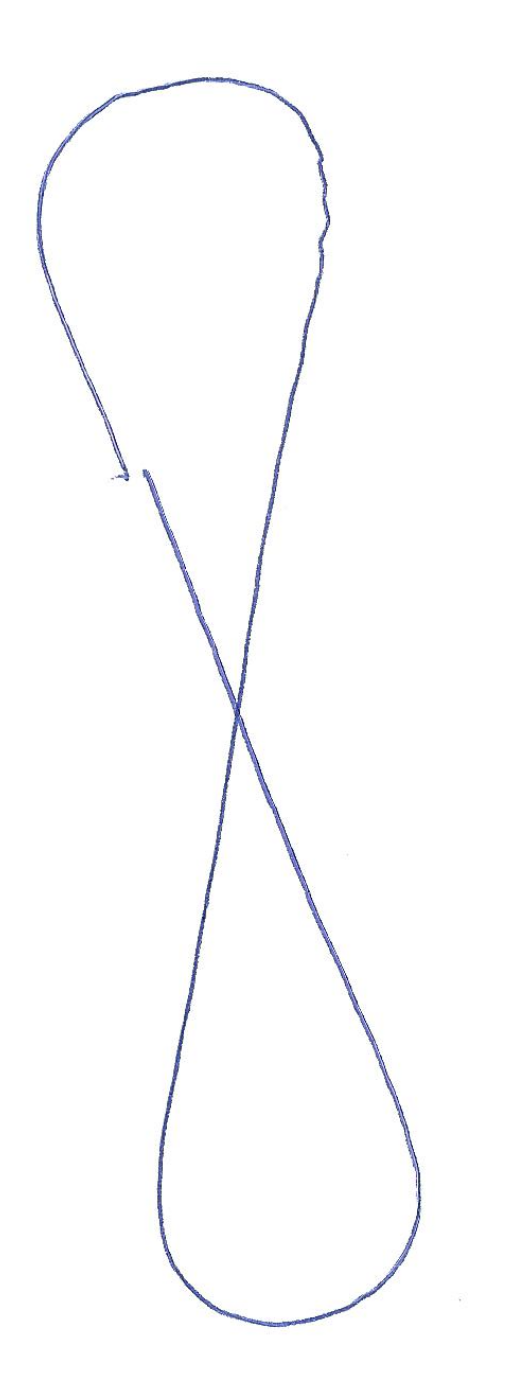}
\caption{Watt's linkage (1784) and the lemniscate it draws.
The pen refill should be inserted in the middle of the white beam.
Note that
four pieces of half-sized beam 6 are used.
}
\label{fig:watt}
\end{center}
\end{figure}

\begin{figure}
\begin{center}
\includegraphics[scale=0.2]{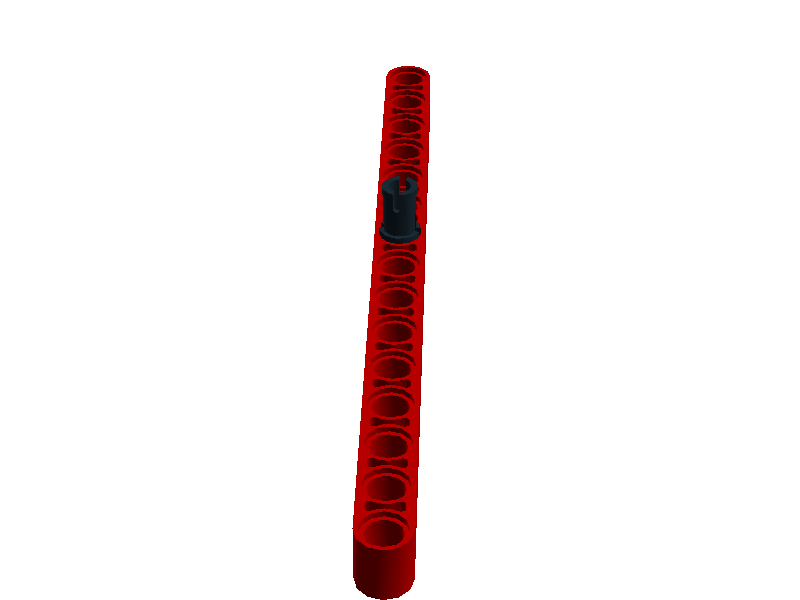}
\includegraphics[scale=0.2]{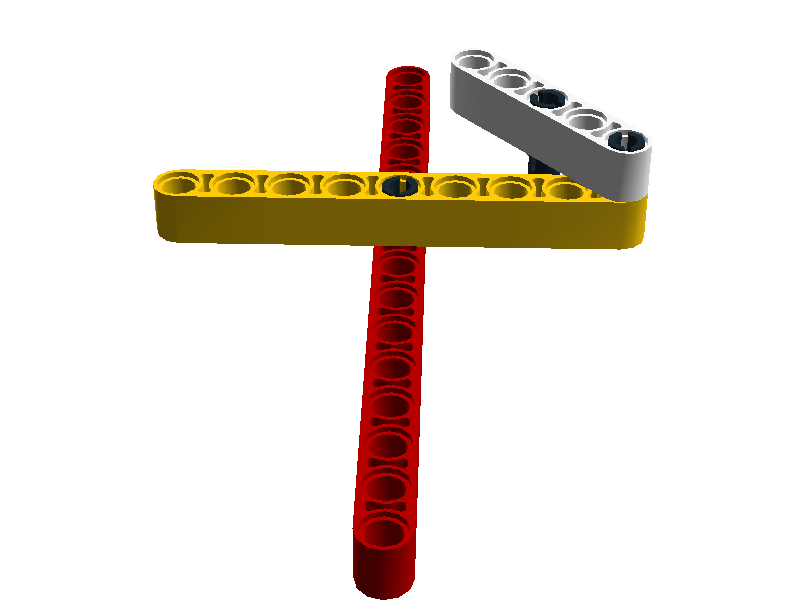} \\
\includegraphics[scale=0.2]{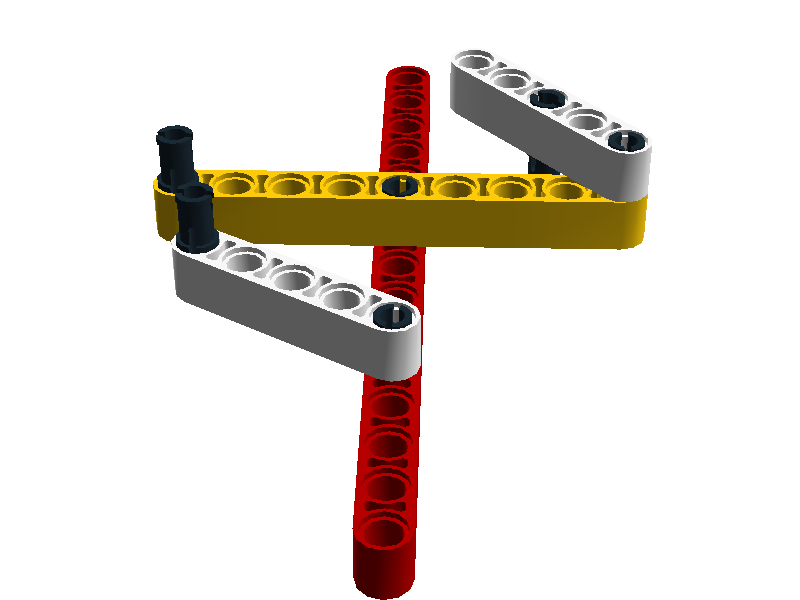}
\includegraphics[scale=0.2]{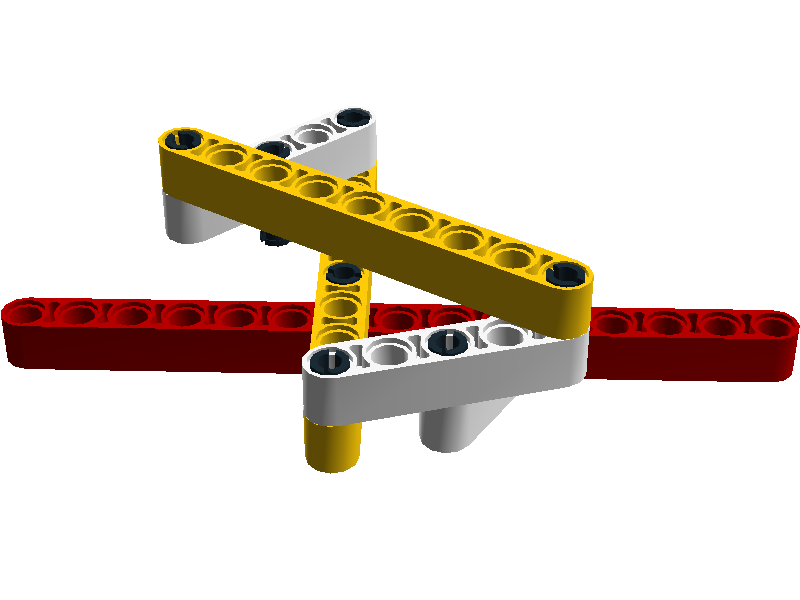}
\caption{Hart's inversor (1874), being constructed in four simple steps.
The pen refill should be inserted in the middle of the white beam at the rear side.}
\label{fig:hart}
\end{center}
\end{figure}

We discuss 5 planar linkages in this paper:
\begin{itemize}
\item a simple ``pair of compasses'' (Fig.~\ref{fig:circle}),
\item Chebyshev's linkage (\Cref{fig:chebyshev,fig:how}),
\item Chebyshev's $\lambda$ mechanism (Fig.~\ref{fig:chebyshevl}),
\item Watt's linkage (Fig.~\ref{fig:watt}) and
\item Hart's inversor (Fig.~\ref{fig:hart}).
\end{itemize}
From the mathematical point of view they are adequately described in the books \cite{kempe,bryantsangwin}.
Here we just give an outline of a possible teaching process by using all of these linkages
to provide an introductory course on planar linkages.

In all models we use LEGO Technic parts (beams and connection pegs) and a pen refill of a specific type
(large capacity, 98.1 millimeters long with a diameter of 6 millimeters, used very commonly in pens).
In all models there is a long red beam with 15 holes (actually its practicable size is just 14 units) which
will be used to firmly hold the model with the one hand (for right-handed pupils this is the left hand),
and a pen refill which will be inserted in a certain connection peg, and hold and moved with the other hand,
on a sheet of paper. This sheet of paper will be pressed down by the whole construction, practically
by the force of the one hand (see Fig.~\ref{fig:how}).

Clearly, the simple compass model (Fig.~\ref{fig:circle}) will draw a circle with the radius of 4 units. This can be introduced
for learners very early, directly after the definition of a circle appears.

The rest of the models require some higher understanding of mathematical concepts. Chebyshev's linkage
(\Cref{fig:chebyshev,fig:how})
is a good example of a ``multi-functional'' model which draws a circle-like curve in its basic form,
but a curve, which looks like a straight line, can also be drawn when the bars are crossed.

The importance of Chebyshev's linkage in the education lies in at least two details. First of all,
it produces two curves which are very similar to well known curves, but they are actually not those (see
\cite{NotLine,bryantsangwin} for a proof). Second, the accuracy of Chebyshev's linkage to these
``similar curves'' is a matter of technical question which is essential in engineering, and ``can mean
the difference between success and failure'' \cite{bryantsangwin}.

For the first impression it is hard to believe that the crossed Chebyshev's linkage does not produce
an exact straight line. For a better understanding it may be useful to build Chebyshev's $\lambda$ mechanism
(Fig.~\ref{fig:chebyshevl})
with the pupils. It produces not just the same curve as Chebyshev's linkage does,
but it can draw the full closed curve in one movement (see
\cite{wiki:chebyshev} for a visual explanation on the equivalence of these two models). By observing the
closed curve, it may be a kind of visual evidence that neither the circle-like part, nor the line-like
part is accurate.

Historically Chebyshev's linkage was not the first one which was identified as a simple machine
producing an almost straight line. Clearly, Watt's linkage (Fig.~\ref{fig:watt})
has been having a much wider attention due to its
general use in several machines, including the steam machine. Despite Watt's linkage does not
produce an exact straight line, either, it is still being used in, for example, suspensions in modern cars
(see \cite{wiki:suspension} for an animated explanation). From the STEM perspective its accuracy
is a ``real life'' technical problem.

Watt's linkage uses almost all parts from our compilation, because of technical reasons.
To ensure a big enough area to move the pen refill freely enough,
we used a helper frame to hold the moving parts. The exact positions of the fixed points of the
yellow components have been ensured by using extra components with 6 holes---they are the hypotenuses
of two Pythagorean triangles to make exact right angles. This model produces about a 7 cm long, seemingly straight
line, but it seems plausible that any part of the resulting curve is not exactly straight because the full
curve is a lemniscate.

The teacher may want to leave enough time for the pupils to do experiments with the
LEGO linkages---maybe without any further instructions.
Small changes in the positions of the connecting holes and pegs may result in various
changes in the output curve. Actually, by doing small modifications on all of the above discussed
linkages, it is impossible to construct a movement having an exact straight line motion.

That is, finally the pupils will \textit{really} want
to build a linkage which can indeed produce an exact straight line. Achieving this is actually
far from being trivial. It can be proven that at least 6 bars are required to make this precisely work,
and to rigorously prove that the model is correct, one may require deeper understanding of some
concepts of planar geometry, including inversions (see \cite{bryantsangwin} for a detailed explanation).
In fact there are several options to build a linkage which produces an exact straight line (see \cite{bryantsangwin}
for other ideas).
In this paper we chose Hart's inversor (Fig.~\ref{fig:hart}) which is moderately difficult to build, but still suitable
for some basic observations.

Unfortunately, this model produces just a very short straight line (on the paper sheet
it is just about 2 centimeters---a longer straight line can be constructed by a model shown
in Fig.~\ref{fig:HartAFrame} in Section \ref{sec:summary}),
and due to some gap between the joining pen refill and the connector peg, the resulting curve
may not look better than the former ones. The only benefit of this fact is to force the pupils to
believe in exact computations and not just in their experience.

\begin{figure}
\begin{center}
\fbox{\includegraphics[width=\textwidth]{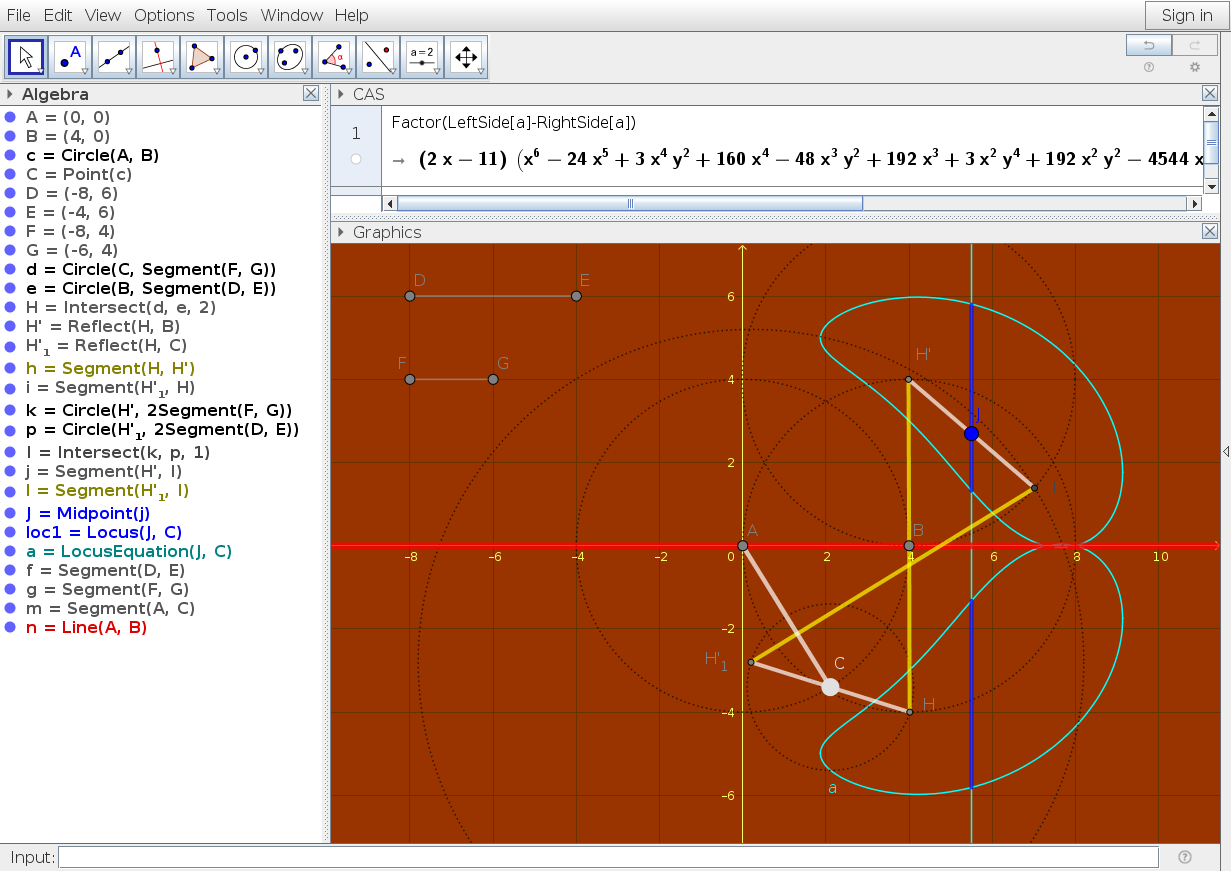}}
\caption{Hart's inversor in GeoGebra. The blue segments, computed numerically, show the expected result of the trace of
the pen refill---there are two disjoint segments because of symmetry.
The cyan curve has been obtained symbolically by using the \texttt{LocusEquation}
command which results in the union of the linear and the sextic component. For mathematical reasons
it is not straightforward to exclude the sextic component automatically: the intersection points of the circles
cannot be explicitly designated by using complex algebraic geometry with the help of Gr\"obner bases.
For similar reasons, a line cannot be restricted to a segment in this approach.
The coloration follows Fig.~\ref{fig:hart}. An online version of this applet
can be found at \cite{Hart-online}---it is possible to change some parameters
by dragging the points $A$, $B$, $D$, $E$, $F$ and $G$, while point $C$ is
constrained to be on the perimeter of circle $c$ (with center $A$ and circumpoint $B$).
}
\label{fig:Hart-gg}
\end{center}
\end{figure}

From the mathematical point of view, Chebyshev's linkage, Chebyshev's $\lambda$ mechanism and also Watt's
linkage is a sextic, that is, an algebraic curve of sixth order. This property is difficult to obtain
by using traditional means in the classroom. Also, it requires some more mathematical background
to understand why a non-linear algebraic curve cannot contain any exactly straight segment parts:
\begin{enumerate}
\item For curves written in the form $y=p(x)=a_0+a_1x+a_2x^2+\ldots+a_nx^n$ this is an
immediate consequence of the fact that the polynomial $p(x)$ has at most $n$ real roots, that is,
two non-zero polynomial functions can have at most $n$ intersections.
\item For arbitrary algebraic curves it is a consequence of B\'ezout's theorem for two dimensions. (See
\cite{CGTA2017} for a short proof.)
\end{enumerate}
The dynamic mathematics software GeoGebra is, however, capable of finding an algebraic equation of
the constructed curves. Since this process in GeoGebra is purely symbolic, the result is reliable,
and no more questions will arise if the construction could produce a straight line or just
an approximation of it. Actually, a technical issue is still important: the equation of the curve
must be converted to a factorized form---this can be however done by using GeoGebra's \texttt{Factor}
command accordingly. When considering Hart's inversor, one can indeed obtain an algebraic
equation of 7th grade, but it splits to a linear and a sextic component (Fig.~\ref{fig:Hart-gg}).

\section{Technical details of the LEGO models}
\label{sec:build}

In this section we give a summary of the used LEGO Technic components. Some details are presented
on the technical challenges on building the models.

Table \ref{tab:set} describes all parts being used in the models. For each part there is a number
on its presence in each model, for example, 8 pieces of pin components are used in Hart's inversor,
but altogether 9 pieces of pins are in the set because 9 ones are needed for Watt's linkage.
Chebyshev's linkage contains 12 parts altogether, but in the set there are 24 different parts
(out of eight kinds of parts: six kinds of beams and two kinds of pins).

In the table there are also prices listed. We show the minimal price for a normal
customer at Brick Owl (\url{https://www.brickowl.com}), a website advertised as a marketplace to
buy and sell LEGO parts, mini-figures and sets. Also, minimal prices of another site, BrickLink (\url{https://www.bricklink.com})
are listed, for customers who want to buy a larger amount of LEGO parts---this simply means that
the lower prices are discounted for wholesale.

\begin{table}
\begin{center}
{\scriptsize
\begin{tabular}{r|cccccccc|c}
Part	& Beam 15 & Beam 11 & Beam 9 & Half beam 6 & Beam 5 & Beam 3 & Long pin\footnote{long pin with friction} &
Pin\footnote{Technic pin with lengthwise friction ridges and center slots}& Total \\
Color	& red & gray\footnote{dark stone gray} & yellow & gray\footnote{dark stone gray} & white &
 gray\footnote{dark stone gray} & blue & black & \\
Code	& 32278	& 32525	& 40490	& 32063	& 32316	& 32523	& 6558	& 2780	& \\
Price\footnote{\label{bo}Brick Owl} &	0.19	& 0.093	& 0.06	& 0.108	& 0.03	& 0.04	& 0.01	& 0.005	& \\
Price\footnote{\label{bl}BrickLink} &	0.1	& 0.02	& 0.03	& 0.04	& 0.0084&0.0084	&0.002	&0.0008	& \\
\hline
Compass	&1	&	&	&	&1	&	&	&1	&3 \\
Cheb.&	1	&2	&1	&	&1	&	&2	&5	&12 \\
Cheb. $\lambda$&	1	&1	&	&2	&	&1	&	&5	&10 \\
Watt	&1	&2	&2	&4	&1	&	&2	&9	&21 \\
Hart	&1	&	&2	&	&3	&	&	&8	&14 \\
\hline
Set	&1	&2	&2	&4	&3	&1	&2	&9	&24 \\
Set price\footnote{Brick Owl}	&0.19	&0.186	&0.12	&0.432	&0.09	&0.04	&0.02	&0.045	&1.123 \\
Set price\footnote{BrickLink}	&0.1	&0.04	&0.06	&0.16	&0.0252	&0.0084	&0.004	&0.0072	& 0.4048\\
\end{tabular}
\caption{Shopping list of the LEGO bricks. Buying components in other colors may be more expensive. Prices
are listed in  \euro .}
\label{tab:set}
} % \small
\end{center}
\end{table}

The table does not show another components of the models, namely the pen refill and the paper sheets.
From these two  clearly the pen refill is the more expensive item, however, various manufacturers have
various products from 0.2 \euro~to more than 3  \euro. Due to our experience the cheaper products are also suitable
for the STEM experiments.

As a total, for a single customer it is possible to buy a complete set of the compilation kit between 1--2 \euro,
without shipping costs.
For classroom use, a teacher may order more parts at once to get complete sets for all pupils at the same time, 
unless they already
have most or all parts from their home sets.

During creation of these models the following restrictions had to be kept in mind:
\begin{itemize}
\item The pen refill has to lie exactly on the paper, and be fixed to a beam hole as much as possible.
For this reason a connector peg is used either at the top or bottom of a hole. This solution is
a home-made way and not officially supported by LEGO, but for the educational purposes it is satisfactory.
\item The red beam (with 15 holes) must lie exactly on the paper, that is, on the ``ground floor''.
Therefore the pen needs to be inserted into a beam being at the height of the first floor (with
a connector peg at the top of the beam, see \Cref{fig:circle,fig:chebyshev})
or the second floor (with a connector peg at the bottom of the beam, see \Cref{fig:chebyshevl,fig:watt,fig:hart}).
\item The ``pair of compasses'' produces an incomplete circle. Of course it is possible to create
a better model, but---to ease educational purposes, that is, to have a general rule how to start a model---the
long red beam is used in this paper to hold the white beam.
Substituting the beams with other ones a perfect circle can be drawn.
\item Some curves can be difficult to draw with one movement. Also Chebyshev's $\lambda$ mechanism may require
some skill to get a continuous curve also at the ``edges''.
\item Since there is no full-sized beam with 6 holes, we need to use two half beams as replacement. Actually,
there exist beams with 7 holes, and they may supplant the beams with 6 holes without using the 7th hole (full-sized
beams are usually cheaper to buy). But, in this case, the model may not look elegant any longer.
\item In many cases some shorter beams can be substituted by longer ones. Sometimes, however, the longer beams
restrict the free movement of other beams.
\item In Watt's linkage it was useful to create exact right angles. We used a Pythagorean triangle with
sides 3, 4 and 5 to achieve this, but there are other methods to mention, see \cite{lego-pyth}.
\item The model of Hart's inversor behaves remarkably when it is moved to its extents. This happens when point $J$
reaches one of the ends of the blue segment (Fig.~\ref{fig:Hart-gg}). At its present state
the model cannot continue its movement, but when removing the connection peg at the pen refill,
the antiparallelogram $HH'IH_1'$ can turn into a parallelogram and point $J$ can continue its way
on the sextic branch (see Fig.~\ref{fig:Hart-gg2}), and a well visible path can be drawn if the
pen refill is inserted again. For mechanical reasons it is not possible to draw the full branch
with this model, however.
\begin{figure}
\begin{center}
\includegraphics[height=5cm]{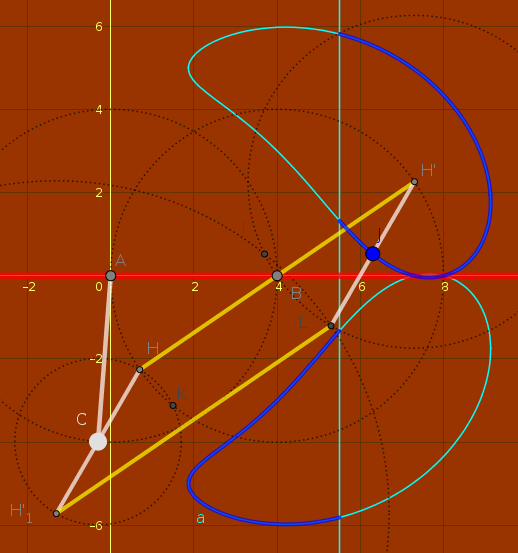}\hskip0.5cm
\includegraphics[trim={0 0.1cm 0 0.3cm},clip,height=5cm]{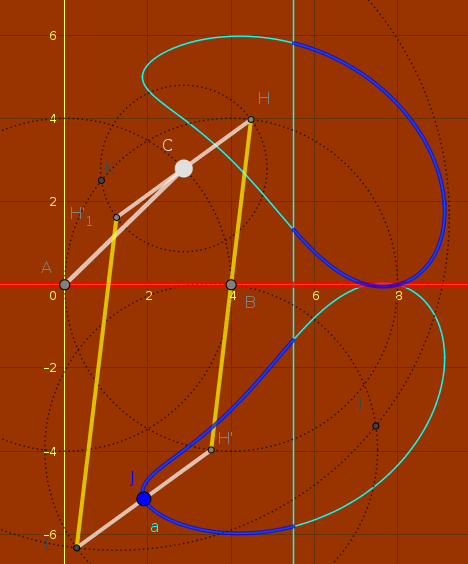}
\caption{The antiparallelogram turns into a parallelogram in Hart's inversor.}
\label{fig:Hart-gg2}
\end{center}
\end{figure}
\item The digital form of the models have been created with the LEGO Digital Designer (LDD) software tool (version 4.3,
available at \url{http://ldd.lego.com}).
It is a very useful way to communicate the steps how a model is built, and finally a step-by-step guide
can also be created. The greatest challenge during creating the electronic models was to align the
connector pegs to the correct holes to link them together in the very final step.
LDD supports this kind of step with its \textit{Hinge} and \textit{Hinge Align} tools.

\end{itemize}

\section{Conclusion, further work}
\label{sec:summary}

We presented a compilation of LEGO Technic parts to promote learning experiments with planar linkages.
A possible didactic approach was sketched up, and some technical details on assembling the parts
were listed.

\begin{figure}
\begin{center}
\includegraphics[height=4.5cm]{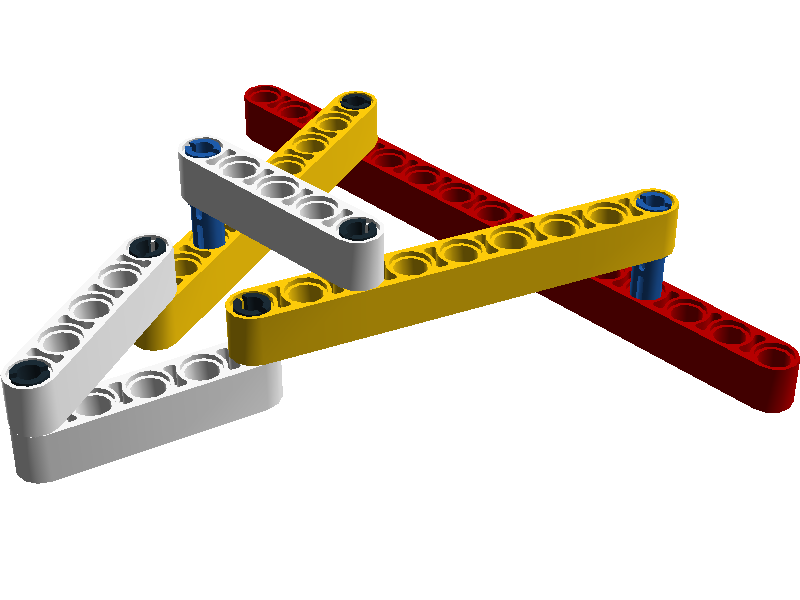}
\includegraphics[height=4.5cm]{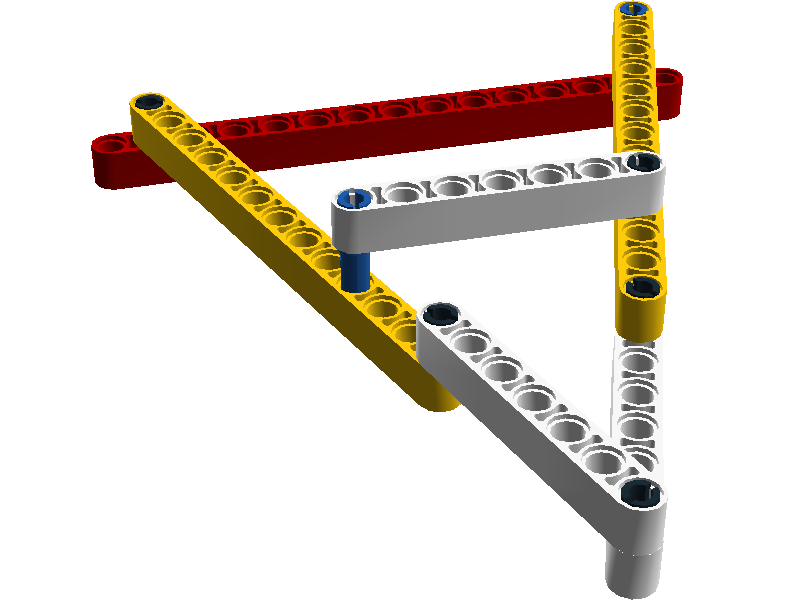}
\caption{Two models to construct Hart's A-frame. The lines that can be drawn are about 4 and 6 centimeters long.
The small one uses only parts from the compilation. The big one contains some longer
beams as well: 13 instead of 9 (the yellow ones) and 7 instead of 5 (the white ones),
not available in our compilation due to save costs. The pen refill has to be inserted at the connection point of the
two white bars.}
\label{fig:HartAFrame}
\end{center}
\end{figure}

The literature of planar linkages is, however, rich in further examples. For instance, several other
types of 4-bar linkages can be built from LEGO parts, and the discussed ones can also be generalized
with different parameters. On the other hand, many other linkages can be built to get an exact
straight line---for further examples we refer to \cite{bryantsangwin}. Two possible examples can be seen
in Fig.~\ref{fig:HartAFrame}---one of them can be built from parts of our compilation as well.

In this paper our approach was to show just one possible way to make the models quickly enough
and available for a wide range of learners. Of course, other construction kits could also be considered to use,
but in our opinion, the LEGO parts are extremely popular and therefore our approach is a very straightforward
way for most learners.

One drawback when using LEGO components is that the parts can only be connected through pins and holes,
and this limits our mathematical model to have only integer lengths appearing---what is more,
some beams (e.g.~a full beam with 6 holes, or more than 15 holes) do not exist as LEGO parts at all,
and lengthening shorter beams may be expensive, unpleasant or not elegant.
This problem could be solved by, for example, 3D printing: the lengths of the bars can be then better fine-tuned.

Modeling 2D and also 3D linkages is an important topic in engineering and more particularly, in robotics.
Today's computers usually do a very good job in sketching up numerically how the artifacts should work
in real life. For high quality products, however, theoretical considerations are also unavoidable:
for example, to decide if a planar curve as the output of a movement really meets the expectations.
In our opinion, using LEGO components, and then extending the collected experience by using computer algebra
and dynamics geometry, can be a fruitful combination in promoting the education process towards
these goals.

\section{Acknowledgments}

The authors are grateful to Chris Sangwin for his useful advice on the topic.
The idea of creating linkages by using GeoGebra and proving some properties on them is suggested by
Tom\'as Recio. The first model of Chebyshev's $\lambda$ mechanism was constructed by Alex Mittermayr.
The figures have been created with the free LDD tool, a product of The LEGO Group.

\bibliography{kovzol,external}

\end{document}